\documentclass[12pt]{article}
\usepackage{amssymb}
\usepackage{mathrsfs}%有更多的数学符号，比如花体$\mathscr P$
\usepackage{pst-poly}     % From pstricks/contrib/pst-poly
\usepackage{cite}

\newtheorem{thm}{Theorem}[section]
\newtheorem{lem}[thm]{Lemma}
\newtheorem{Def}[thm]{Definition}
\newtheorem{cor}[thm]{Corollary}

\newtheorem{rem}[thm]{Remark}

\newenvironment{pf}[1][Proof]{\noindent\textbf{#1.} }{\hfill\rule{1mm}{2mm}}

\makeatletter \@addtoreset{equation}{section} \makeatother

\begin{document}
\title{\bf Generalized Measures of Fault
Tolerance in $(n,k)$-star Graphs\thanks {The work was supported by
NNSF of China (No.11071233).}}
\author
{Xiang-Jun Li\quad
%$^a$\quad Yun-Chao Wei$^{a,b}$\footnote{When the joint
%work was completed the author was visiting School of Mathematical
%Sciences, University of Science and Technology of China.}\quad
Jun-Ming
Xu\footnote{Corresponding author: xujm@ustc.edu.cn (J.-M. Xu)}\\
$^a${\small School of Mathematical Sciences, University of Science
and Technology of China,}\\
{\small  Wentsun Wu Key Laboratory of CAS, Hefei, 230026, China}  \\
%$^b${\small College of Information Technology, Shanghai Ocean
%University, Shanghai, 201306, China}
 }
\date{}
 \maketitle

%\begin{frontmatter}
\begin{abstract}

This paper considers a kind of generalized measure $\kappa_s^{(h)}$
of fault tolerance in the $(n,k)$-star graph $S_{n,k}$ and
determines $\kappa_s^{(h)}(S_{n,k})=n+h(k-2)-1$ for $2 \leqslant k
\leqslant n-1$ and $0\leqslant h \leqslant n-k$, which implies that
at least $n+h(k-2)-1$ vertices of $S_{n,k}$ have to remove to get a
disconnected graph that contains no vertices of degree less than
$h$. This result contains some known results such as Yang et al.
[Information Processing Letters, 110 (2010), 1007-1011].

\vskip6pt


\end{abstract}

\section{Introduction}
It is well known that interconnection networks play an important
role in parallel computing/communication systems. An interconnection
network can be modeled by a graph $G=(V, E)$, where $V$ is the set
of processors and $E$ is the set of communication links in the
network.
%For all the graph terminologies and notations not defined
%here, we follow\cite{x03}. Then we use graphs and networks
%interchangeably in this paper.
The connectivity $\kappa(G)$ of a graph $G$ is an important
measurement for fault-tolerance of the network, and the larger
$\kappa(G)$ is, the more reliable the network is.

A subset of vertices $S$ of a connected graph $G$ is called a {\it
vertex-cut} if $G-S$ is disconnected. The {\it connectivity}
$\kappa(G)$ of $G$ is defined as the minimum cardinality over all
vertex-cuts of $G$. Because $\kappa$ has many shortcomings, one
proposes the concept of the $h$-super connectivity of $G$, which can
measure fault tolerance of an interconnection network more
accurately than the classical connectivity $\kappa$.

A subset of vertices $S$ of a connected graph $G$ is called an {\it
$h$-super vertex-cut}, or {\it $h$-cut} for short, if $G-S$ is
disconnected and has the minimum degree at least $h$. The {\it
$h$-super connectivity} of $G$, denoted by $\kappa^{(h)}_s(G)$, is
defined as the minimum cardinality over all $h$-cuts of $G$. It is
clear that, if $\kappa_s^{(h)}(G)$ exists, then
 $$
 \kappa(G)=\kappa_s^{(0)}(G)\leqslant \kappa_s^{(1)}(G)\leqslant \kappa_s^{(2)}(G)\leqslant
 \cdots \leqslant  \kappa_s^{(h-1)}(G)\leqslant \kappa_s^{(h)}(G).
 $$

For any graph $G$ and integer $h$, determining $\kappa_s^{(h)}(G)$
is quite difficult. In fact, the existence of $\kappa_s^{(h)}(G)$ is
an open problem so far when $h\geqslant 1$. Only a little knowledge
of results have been known on $\kappa_s^{(h)}$ for particular
classes of graphs and small $h$'s.

This paper is concerned about $\kappa_s^{(h)}$ for the $(n,k)$-star
graph $S_{n,k}$. For $k=n-1$, $S_{n,n-1}$ is isomorphic to a star
graph $S_n$, Cheng and Lipman~\cite{cl02}, Hu and Yang~\cite{hy97},
Nie {\it et al.}~\cite{nlx04} and Rouskov {\it et al.}~\cite{rls96},
independently, determined $\kappa_s^{(1)}(S_n)=2n-4$ for $n\geqslant
3$. Very recently, Yang {\it et al.}~\cite{ylg10} have showed that
if $2\leqslant k\leqslant n-2$ then $\kappa_s^{(1)}(S_{n,k})=n+k-3$
for $n\geqslant 3$ and $\kappa_s^{(2)}(S_{n,k})=n+2k-5$ for
$n\geqslant 4$.

We, in this paper, will generalize these results by proving that
$\kappa _s^{(h)}(S_{n,k})=n+h(k-2)-1$ for $2\leqslant k\leqslant
n-1$ and $0\leqslant h\leqslant n-k$.

The proof of this result is in Section 3. In Section 2, we recall
the structure of $S_{n,k}$ and some lemmas used in our proofs.

\section{Definitions and lemmas}

For given integer $n$ and $k$ with $1\leqslant k\leqslant n-1$, let
$I_n=\{1,2,\ldots,n\}$ and $P(n,k)=\{ p_{1}p_{2}\ldots p_{k}:\
p_{i}\in I_n, p_{i}\neq p_{j}, 1\leqslant i\neq j\leqslant k\}$, the
set of $k$-permutations on $I_n$. Clearly, $|P(n,k)|=n\,!/(n-k)\,!$.

\begin{Def}\label{def2.1}
The $(n,k)$-star graph $S_{n,k}$ is a graph with vertex-set
$P(n,k)$. The adjacency is defined as follows: a vertex
$p=p_{1}p_{2}\ldots p_{i}\ldots p_{k}$ is adjacent to a vertex

(a)\ $p_{i}p_{2}\cdots p_{i-1}p_{1}p_{i+1}\cdots p_{k}$, where
$2\leqslant i\leqslant k$ (swap $p_{1}$ with $p_{i}$).

(b)\ $\alpha p_{2}p_{3}\cdots p_{k}$, where $\alpha\in I_n\setminus
\{p_{i}:\ 1\leqslant i\leqslant k\}$ (replace $p_{1}$ by $\alpha$).
\end{Def}

The vertices of type $(a)$ are referred to as {\it swap-neighbors}
of $p$ and the edges between them are referred to as {\it swap-edge}
or {\it $i$-edges}. The vertices of type $(b)$ are referred to as
{\it unswap-neighbors} of $p$ and the edges between them are
referred to as {\it unswap-edges}. Clearly, every vertex in
$S_{n,k}$ has $k-1$ swap-neighbors and $n-k$ unswap-neighbors.
Usually, if $x=p_1p_2\dots p_k$ is a vertex in $S_{n,k}$, we call
$p_i$ the {\it $i$-th bit} for each $i\in I_k$.

The $(n,k)$-star graph $S_{n,k}$ is proposed by Chiang and
Chen~\cite{cc95} who showed that $S_{n,k}$ is $(n-1)$-regular
$(n-1)$-connected.

\begin{lem}\label{lem2.2}
For any $\alpha=p_1p_2\cdots p_{k-1}\in P(n,k-1)$ $(k \geqslant 2)$,
let $V_\alpha=\{p\alpha:\ p\in I_n\setminus \{p_i:\ i\in
I_{k-1}\}\}$. Then the subgraph of $S_{n,k}$ induced by $V_\alpha$
is a complete graph of order $n-k+1$, denoted by $K^\alpha_{n-k+1}$.
\end{lem}

\begin{pf}
For any two vertices $p\alpha$ and $q\alpha$ in $V_\alpha$ with
$p\ne q$, by the condition $(b)$ of Definition~\ref{def2.1},
$p\alpha$ and $q\alpha$ are linked in $S_{n,k}$ by an unswap-edge.
Thus, the subgraph of $S_{n,k}$ induced by $V_\alpha$ is a complete
graph $K_{n-k+1}$.
\end{pf}

\vskip6pt

By Lemma~\ref{lem2.2}, the vertex-set $P(n,k)$ of $S_{n,k}$ can be
decomposed into $|P(n,k-1)|$ subsets, each of which induces a
complete graph $K_{n-k+1}$. It is clear that, for any two distinct
elements $x$ and $y$ in $P(n,k)$, if they are in different complete
subgraphs $K^\alpha_{n-k+1}$ and $K^\beta_{n-k+1}$ $(\alpha\ne
\beta)$, then there is at most %one swap-edge
one edge between $x$ and $y$ in
$S_{n,k}$, which is %an $i$-edge if and only if $\alpha$ and $\beta$
%differ in only the $i$-th bit.
a swap-edge if and only if $\alpha$ and $\beta$
differ in only one bit. Thus, we have the following
conclusion.

\begin{lem}\label{lem2.3}
The vertex-set of $S_{n,k}$ can be partitioned into $|P(n,k-1)|$
subsets, each of which induces a complete graph of order $n-k+1$.
Furthermore, there is at most one swap edge between any two complete
graphs.
\end{lem}

Let $S^{t:i}_{n-1,k-1}$ denote a subgraph of $S_{n,k}$ induced by
all vertices with the $t$-th bit $i$ for $2\leqslant t\leqslant k$.
The following lemma is a slight modification of the result of Chiang
and Chen~\cite{cc95}.

\begin{lem}\label{lem2.4}%\textnormal{(Chiang and Chen~\cite{cc95})}
For a fixed integer $t$ with $2\leqslant t\leqslant k$, $S_{n,k}$
can be decomposed into $n$ subgraphs $S^{t:i}_{n-1,k-1}$, which is
isomorphic to $S_{n-1,k-1}$, for each $i\in I_n$. Moreover, there
are $\frac{(n - 2)!}{(n - k)!}$ independent swap-edges between
$S^{t:i}_{n-1,k-1}$ and $S^{t:j}_{n-1,k-1}$ for any $i,j\in I_n$
with $i\ne j$.
\end{lem}

%\begin{pf}
%If we remove the $t$-th bits of all vertices in $S^{t:i}_{n-1,k-1}$,
%we obtain an $S_{n-1,k-1}$. That is, $S^{t:j}_{n-1,k-1}$ is
%isomorphic to $S_{n-1,k-1}$. So, we only need to show
%$S^{t:i}_{n-1,k-1}, i\not =j$ and $S^{t:l}_{n-1,k-1}$ are
%isomorphic. We define the one-to-one mapping $\rho$ in $I_n$:

%$\rho(i) = l, p(l) = i,$

%$\rho(x) =x$ for $x\in  I_n - \{i,l\}.$

%We also define a bijection $ M$ by:

%$M(p_1p_2\ldots p_k) = \rho(p_1)\rho(p_2)\ldots\rho(p_k)$ for a
%vertex $p_1p_2\ldots p_k$ in $S_{n,k}$

%Obviously, $M$ transforms the vertices of $S^{t:i}_{n-1,k-1}$, into
%those of $S^{t:l}_{n-1,k-1}$, and preserves adjacency. It is easy to
%check  there are $\frac{(n - 2)!}{(n - k)!}$ independent swap-edges
%between $S^{t:i_1}_{n-1,k-1}$ and $S^{t:i_2}_{n-1,k-1}$ for any
%$i_1,i_2\in I_n$ with $i_1\ne i_2$.
%\end{pf}

\begin{lem}\label{lem2.6}
\textnormal{(Chen {\it et al.}~\cite{cdyf08})} In an $S_{n,k}$, a
cycle has a length at least 6 if it contains a swap-edge.
\end{lem}

%\newpage

\section{Main results}

In this section, we present our main results, that is, we determine
the $h$-super connectivity of the $(n,k)$-star graph $S_{n,k}$.
Since $S_{n,1}\cong K_n$, we only consider the case of $k\geqslant
2$ in the following discussion.

\begin{lem}\label{lem3.1}
$\kappa_s^{(h)}(S_{n,k})\leqslant n+h(k-2)-1$ for $2 \leqslant k \le
n-1$ and $0\leqslant h \leqslant n-k$.
\end{lem}

\begin{pf}
By our hypothesis of $h \leqslant n-k$, for any $\alpha\in
P(n,k-1)$, we can choose a subset $X\subseteq V(K^\alpha_{n-k+1})$
such that $|X|=h+1$. Then the subgraph of $K^\alpha_{n-k+1}$ induced
by $X$ is a complete graph $K_{h+1}$. Let $S$ be the neighbor-set of
$X$ in $S_{n,k}-X$. Clearly, $V(K^\alpha_{n-k+1}-X)\subseteq S$,
that is, $X$ has exactly $n-k+1-|X|$ unswap-neighbors in
$V(K^\alpha_{n-k+1}-X)\cap S$. Since $S_{n,k}$ is $(n-1)$-regular,
every vertex of $X$ has exactly $(k-1)$ swap-neighbors are not in
$K^\alpha_{n-k+1}$. Moreover, any two swap-neighbors of $X$ are
different from each other by Lemma~\ref{lem2.3}. It follows that
 \begin{equation}\label{e3.1}
 |S|=n-k+1-|X|+|X|(k-1)=n+h(k-2)-1.
 \end{equation}

We now need to show that $S$ is an $h$-cut of $S_{n,k}$. Clearly,
$S$ is a vertex-cut of $S_{n,k}$ since $S_{n,k}$ is not a complete
graph for $k\geqslant 2$. We only need to show that every vertex of
$S_{n,k}-(X\cup S)$ has degree at least $h$. Let $u$ be a vertex in
$S_{n,k}-(X\cup S)$. If $u$ has a neighbor $v$ in $S\cap
V(K^\alpha_{n-k+1})$, then $u$ is a swap-neighbor of $v$ since all
the unswap-neighbors of $v$ are in $V(K^\alpha_{n-k+1})$. If $u$ has a neighbor $v$ in
$S\setminus V(K^\alpha_{n-k+1})$, then $v$ has a swap-neighbor in
$V(K^\alpha_{n-k+1})$. Moreover, if $u$ has two neighbor $v,v'$ in
$S$, then three vertices $u, v$ and $v'$ are concluded in a cycle of
length at most $5$ and containing at least one swap-edge, which
contradicts with Lemma~\ref{lem2.6}. Thus, $u$ has at most one
neighbor in $S$. In other words, $u$ has at least $n-2$ neighbor in
$S_{n,k}-S$. Since $n-2\geqslant n-k\geqslant h$ for $k\geqslant 2$,
%and $n-k \geqslant h$,
$u$ has degree at least $h$ in $S_{n,k}-S$.
By the arbitrariness of $u\in S_{n,k}-(X\cup S)$, $S$ is an $h$-cut
of $S_{n,k}$, and so
 $$
 \kappa _s^{(h)}(S_{n,k})\leqslant |S|=n+h(k-2)-1
 $$
as required. The lemma follows.
\end{pf}

\begin{cor}\label{cor3.2}
$\kappa_s^{(h)}(S_{n,2})=n-1$ for $0 \leqslant h\leqslant n-2$.
\end{cor}

\begin{pf}\label{cor3.1}
On the one hand, $\kappa _s^{(h)} (S_{n,2})\leqslant n-1$ by
Lemma~\ref{lem3.1} when $k=2$. On the other hand, $\kappa _s^{(h)}
(S_{n,2})\ge \kappa(S_{n,2})=n-1$.
\end{pf}

\vskip6pt

To state and prove our main results, we need some notations. Let $S$
be an $h$-cut of $S_{n,k}$ and $X$ be the vertex-set of a connected
component of $S_{n,k}-S$.
%Let $t\in I_n-\{1\}$ such
%that $|\{i\in I_n: X\cap S^{t:i}_{n-1,k-1}\not=\emptyset \}| $ as
%much as possible,  and denote $S^{t:i}_{n-1,k-1}$ simply as
%$S^i_{n-1,k-1}$ for $i\in I_n$.
For a fixed $t\in I_k\setminus\{1\}$ and any $i\in I_n$, let
 \begin{equation}\label{e3.2}
 \begin{array}{l}
 Y=V(S_{n,k}-S-X),\\
 X_i=X\cap V(S^{t:i}_{n-1,k-1}),\\
 Y_i=Y\cap V(S^{t:i}_{n-1,k-1})\ {\rm and}\\
 S_i=S\cap V(S^{t:i}_{n-1,k-1}),
 \end{array}
\end{equation}
and let
\begin{equation}\label{e3.3}
\begin{array}{l}
 J=\{i\in I_n:\ X_i\ne\emptyset\},\\
 J'=\{i\in  J:\ Y_i\not=\emptyset\}\ \ {\rm and} \\
 T=\{i\in I_n:\  Y_i\ne\emptyset\}.
 \end{array}
 \end{equation}

\begin{lem}\label{lem3.3}
Let $S$ be a minimum $h$-cut of $S_{n,k}$ and $X$ be the vertex-set
of a connected component of $S_{n,k}-S$.
%Let $t\in I_n\setminus\{1\}$ such that $|J|$ is as large as possible.
%, and denote $S^{t:i}_{n-1,k-1}$ simply as
%$S^i_{n-1,k-1}$ for $i\in I_n$. For each $i\in I_n$, let $X_i=X\cap
%V(S^i_{n-1,k-1}), Y_i=Y\cap V(S^i_{n-1,k-1}),  S_i=S\cap
%V(S^i_{n-1,k-1}), J=\{i\in I_n:\ X_i\ne\emptyset\}$, $J'=\{i\in
%J:X_i\not=\emptyset \not= Y_i\}$ and $ T=\{i\in I_n:\
%Y_i\ne\emptyset\}$.
If $3\leqslant k\leqslant n-1$ and $1\leqslant h\leqslant n-k$ then,
for any $t\in I_k\setminus\{1\}$,

{\rm (a)}\ $S_i$ is an $(h-1)$-cut of $S^{t:i}_{n-1,k-1}$ for any
$i\in J'$,

{\rm (b)}\ $\kappa_s^{(h)}(S_{n,k})\geqslant |J'|\
\kappa_s^{(h-1)}(S_{n-1,k-1})$,

{\rm (c)}\ $J\cup T =I_n$.

\end{lem}

\begin{pf}
(a)\ By the definition of $J'$, $S_i$ is a vertex-cut of
$S^{t:i}_{n-1,k-1}$ for any $i\in J'$. For any vertex $x$ in
$S^{t:i}_{n-1,k-1}-S_i$, since $x$ has degree at least $h$ in
$S_{n,k}-S$ and has exactly one neighbor outsider
$S^{t:i}_{n-1,k-1}$, $x$ has degree at least $h-1$ in
$S^{t:i}_{n,k}-S_i$. This fact shows that $S_i$ is an $(h-1)$-cut of
$S^{t:i}_{n-1,k-1}$ for any $i\in J'$.

(b)\ By the assertion (a), we have
$|S_i|\geqslant\kappa_s^{(h-1)}(S_{n-1,k-1})$, and so
 $$
 \kappa_s^{(h)}(S_{n,k})=|S|\geqslant \sum_{i\in J'} |S_i|\geqslant
 |J'|\kappa_s^{(h-1)}(S_{n-1,k-1}).
  $$

(c)\ If $J\cup T \not=I_n$, that is, $I_n\setminus(J\cup T)\not =
\emptyset$, then there exists an $i_0\in I_n$ such that
$V(S^{t:i_0}_{n-1,k-1})=S_{i_0}$. Thus, we have
 $$
\begin{array}{rl}
\kappa_s^{(h)}(S_{n,k}) &=|S|\geqslant |S_{i_0}| = \frac{(n-1)!}{(n-k)!}\\
 &\geqslant (n-1)(n-2)\\
 &=n+(n-1)(n-3)-1\\
 &>n+(n-3)(n-3)-1\\
 &\geqslant n+h(k-2)-1,
   \end{array}
 $$
which contradicts to Lemma~\ref{lem3.1}. Thus, $J\cup T =I_n$. The
Lemma follows.
\end{pf}

\begin{thm}\label{thm3.4}
$\kappa _s^{(h)}(S_{n,k})=n+h(k-2)-1$ for $2\leqslant k\leqslant
n-1$ and $0\leqslant h\leqslant n-k$.
\end{thm}

\begin{pf}
By Lemma~\ref{lem3.1}, we only need to prove that, for $2\leqslant
k\leqslant n-1$ and $0\leqslant h\leqslant n-k$,
 \begin{equation}\label{e3.4}
\kappa _s^{(h)}(S_{n,k})\geqslant n+h(k-2)-1.
 \end{equation}

We proceed by induction on $k\geqslant 2$ and $h\geqslant 0$. The
inequality (\ref{e3.4}) is true for $k=2$ and any $h$ with $0
\leqslant h\leqslant n-2$ by Corollary~\ref{cor3.2}. The inequality
(\ref{e3.4}) is also true for $h=0$ and any $k$ with $2\leqslant
k\leqslant n-1$ since $\kappa
_s^{(0)}(S_{n,k})=\kappa(S_{n,k})=n-1$. Assume the induction
hypothesis for $k-1$ with $k\geqslant 3$ and for $h-1$ with
$h\geqslant 1$, that is,
 %\begin{equation}\label{e3.5}
%\kappa _s^{(h)}(S_{n-1,k-1})\geqslant n+h(k-3)-2,
 %\end{equation}
%and
 \begin{equation}\label{e3.6}
\kappa _s^{(h-1)}(S_{n-1,k-1})\geqslant n+(h-1)(k-3)-2.
 \end{equation}

Let $S$ be a minimum $h$-cut of $S_{n,k}$ and $X$ be the vertex-set
of a minimum connected component of $S_{n,k}-S$. Use notations
defined in (\ref{e3.2}) and (\ref{e3.3}). Choose $t\in
I_k\setminus\{1\}$ such that $|J|$ is as large as possible. For each
$i\in I_n$, we write $S^i_{n-1,k-1}$ for $S^{t:i}_{n-1,k-1}$ for
short. We consider three cases depending on $|J'|=0$, $|J'|=1$ or
$|J'|\geqslant2$.

\vskip6pt

{\bf Case 1.}\  $|J'|=0$,

In this case, $X_i\ne\emptyset$ and $Y_i=\emptyset$ for each $i\in
J$, that is, $J\cap T=\emptyset$. By Lemma~\ref{lem3.3} (c),
$|J|\geqslant 2$ or $|T|\geqslant 2$ since $n\geqslant 4$. Clearly,
$J\ne\emptyset$ and $T\ne\emptyset$. Without loss of generality,
assume $|J|\geqslant 2$, $\{i_{1},i_2\}\subseteq J$ and $i_3\in T$.
By Lemma~\ref{lem2.4}, there are $\frac{(n - 2)!}{(n - k)!}$
independent swap-edges between $S^{i_1}_{n-1,k-1}$ (resp.
$S^{i_2}_{n-1,k-1}$) and $S^{i_3}_{n-1,k-1}$, each edge of which has
at least one end-vertex in $S$. Since $J\cap T=\emptyset$ and
$S_{i_1}\cap S_{i_2}=\emptyset$, we have that
\begin{equation}\label{e3.7}
\begin{array}{c}
 |S|\geqslant 2\ \frac{(n-2)!}{(n-k)!}.
  \end{array}
 \end{equation}
Noting that, for $k=3$,
 $$
 \begin{array}{c}
2\ \frac{(n-2)!}{(n-k)!}\geqslant 2(n-2)\geqslant n+(n-3)-1\geqslant
n+h(k-2)-1,
  \end{array}
 $$
and, for $k\geqslant 4$,
 $$
 \begin{array}{c}
 2\ \frac{(n-2)!}{(n-k)!}\geqslant 2(n-2)(n-3)\geqslant n+(n-3)(n-3)-1\geqslant
 n+h(k-2)-1,
  \end{array}
 $$
we have that
\begin{equation}\label{e3.8}
 \begin{array}{c}
 2\ \frac{(n-2)!}{(n-k)!}\geqslant
 n+h(k-2)-1\ \ {\rm for}\ k\geqslant 3.
   \end{array}
 \end{equation}
It follows from (\ref{e3.7}) and (\ref{e3.8}) that
$$
\begin{array}{rl}
\kappa_s^{(h)}(S_{n,k})& =|S|\geqslant 2\ \frac{(n-2)!}{(n-k)!}
 \geqslant n+h(k-2)-1.
 \end{array}
$$

\vskip6pt

{\bf Case 2.}\  $|J'|=1$,

Without loss of generality, assume %$Y_1\not =\emptyset$.
$J'=\{1\}$. By
Lemma~\ref{lem3.3} (a), $S_1$ is an $(h-1)$-cut of $S^1_{n-1,k-1}$.
Let $S'=S\setminus S_1$.

If $|S'|\geqslant n-2$ then, by (\ref{e3.6}), we have that
$$
\begin{array}{rl}
\kappa_s^{(h)}(S_{n,k})& =|S| = |S_1|+|S'| \geqslant \kappa_s^{(h-1)}(S_{n-1,k-1})+(n-2) \\
&\geqslant (n+(h-1)(k-3)-2)+(n-2)\\
&\geqslant (n+(h-1)(k-3)-2)+(h+k-2)\\
&= n+h(k-2)-1
\end{array}
    $$

We now assume $|S'|\leqslant n-3$. We claim $|J|=1$. Suppose to the
contrary $|J|\geqslant 2$.

If $|T|= 1$ then, by Lemma 3.3 (c), we have $|J|=n$. Then
$|S_i|\geqslant 1$ for $i\in J$, otherwise there exists $i\in
J\setminus J'$ such that $X_i=V(S^i_{n-1,k-1})$, then $|X| >
|X_i|=|V(S^i_{n-1,k-1})|>|Y|$, which contradicts to  the minimality
of $X$. Therefore,  $|S'|\geqslant n-1$, a contradiction.

If $|T|\geqslant 2$, assume that $i_1\in J\setminus J'$ and $i_2\in
T\setminus J'$, then $X_{i_1}\not=\emptyset,
Y_{i_1}=\emptyset,X_{i_2}=\emptyset, Y_{i_2}\not=\emptyset$. By
Lemma~\ref{lem2.4}, there are $\frac{(n - 2)!}{(n - k)!}$
independent swap-edges between $S^{i_1}_{n-1,k-1}$ and
$S^{i_2}_{n-1,k-1}$, each edge of which must have one end-vertex in
$S'$. Thus, we have
 $$%\begin{equation}\label{e3.7}
 \begin{array}{c}
|S'|\geqslant \frac{(n- 2)!}{(n - k)!}\geqslant n-2\ \ {\rm for}\
k\geqslant 3,
 \end{array}
 $$%\end{equation}
a contradiction.

Thus, $|J|=1$. %Without loss of generality, assume $J=\{1\}$.
We have $J=\{1\}$ since $\{1\}=J'\subseteq J$. Then
$X_1=X$ and $|X_1|\geqslant h+1$. By the choice of $t$, the $i$-th
($i\ne 1$) bits of all vertices in $X_1$ are same, and so $X_1$ a
complete graph. Thus, as computed in (\ref{e3.1}), we have that
     $$
\begin{array}{rl}
\kappa_s^{(h)}(S_{n,k})& =|S|=n+(|X_1|-1)(k-2)-1\geqslant
n+h(k-2)-1.
\end{array}
    $$

\vskip6pt

{\bf Case 3.}\  $|J'|\geqslant 2$.

By Lemma 3.3 (b) and (\ref{e3.6}), we have that
$$
\begin{array}{rl}
\kappa_s^{(h)}(S_{n,k})=|S| &\geqslant |J'|\kappa_s^{(h-1)}(S_{n-1,k-1})\\
 &\geqslant 2(n+(h-1)(k-3)-2)\\
 &\geqslant n+(h+k)+2(h-1)(k-3)-4\\
 &=n+h(k-2)+(h-1)(k-3)-1\\
 &\geqslant n+h(k-2)-1.
 \end{array}
 $$

By the induction principle, the theorem follows.
\end{pf}

\begin{cor}\textnormal{(Yang et al.~\cite{ylg10})}
If $2\leqslant k\leqslant n-2$ then $\kappa_s^{(1)}(S_{n,k})=n+k-3$
for $n\geqslant 3$ and $\kappa_s^{(2)}(S_{n,k})=n+2k-5$ for
$n\geqslant 4$.
\end{cor}

As we have known, when $k=n-1$, $S_{n,n-1}$ is isomorphic to the
star graph $S_n$. Akers and Krishnamurthy~\cite{ak89} determined
$\kappa(S_n)=n-1$ for $n\geqslant 2$; Cheng and Lipman~\cite{cl02}, Hu and Yang~\cite{hy97},
Nie {\it et al.}~\cite{nlx04} and Rouskov {\it et al.}~\cite{rls96},
independently, determined
$\kappa_s^{(1)}(S_n)=2n-4$ for $n\geqslant 3$. All these results can
be obtained from our result by setting $k=n-1$ and $h=0,1$,
respectively.

\begin{cor}
$\kappa(S_n)=n-1$ for $n\geqslant 2$ and $\kappa_s^{(1)}(S_n)=2n-4$
for $n\geqslant 3$.
\end{cor}

\begin{rem}
Wan and Zhang~\cite{wz09} determined $\kappa_s^{(2)}(S_n)=6(n-3)$
for $n\geqslant 4$. Thus, our result is invalid for $\kappa_s^{(h)}(S_n)$
when $h\geqslant 2$. Thus, determining $\kappa_s^{(h)}(S_n)$ for
$h\geqslant 2$ needs other technique.
\end{rem}

\end{document}